\documentclass[a4paper]{article}

\usepackage{amsmath, amsfonts, amsthm}
\usepackage{graphicx}

\newtheorem{thm}{Theorem}[section]

\newtheorem{cor}[thm]{Corollary}
\newtheorem{bijection}{Bijection}

\theoremstyle{definition}
\newtheorem{exm}[thm]{Example}

\newcommand{\Pos}{\mathbb{P}}
\newcommand{\N}{\mathbb{N}}
\newcommand{\Z}{\mathbb{Z}}
\newcommand{\comp}{\mathcal{C}}

\newcommand{\abs}[1]{\left\lvert #1 \right\rvert}
\newcommand{\rev}[1]{\left\langle #1 \right\rangle}

\newcommand{\mn}{\mbox{-}}

\title{ %
    Enumerating Segmented Patterns in Compositions and Encoding by
    Restricted Permutations
}

\author{
    Sergey Kitaev\thanks{Department of Mathematics, University of
    California, San Diego, CA 92093, USA. Supported by the Institut
    Mittag--Leffler, Sweden, in Winter 2005.} \\
    \and Tyrrell B. McAllister\thanks{Mathematics, University of
    California, Davis, CA 95616, USA. Supported by the NSF while
    staying at the Institut Mittag--Leffler, Sweden, in Winter 2005.}
    \and T. Kyle Petersen\thanks{Department of Mathematics, Brandeis University, Waltham, MA 02454, USA. Supported by the NSF while staying at
    the Institut Mittag--Leffler, Sweden, in Winter 2005.}
}

\begin{document}
    \maketitle

    \begin{abstract}
        \noindent A {\em composition} of a nonnegative integer \( n \)
        is a sequence of positive integers whose sum is \( n \).  A
        composition is {\em palindromic} if it is unchanged when its
        terms are read in reverse order.  We provide a generating
        function for the number of occurrences of arbitrary segmented
        partially ordered patterns among compositions of \( n \) with
        a prescribed number of parts.  These patterns generalize the
        notions of rises, drops, and levels studied in the literature.
        We also obtain results enumerating parts with given sizes and
        locations among compositions and palindromic compositions with
        a given number of parts.  Our results are motivated by
        ``encoding by restricted permutations," a relatively
        undeveloped method that provides a language for describing
        many combinatorial objects.  We conclude with some examples
        demonstrating bijections between restricted permutations and
        other objects.
    \end{abstract}

    \section{Introduction} \label{intro}

    A {\em composition} of a nonnegative integer $n$ is a sequence \(
    \alpha = \alpha_{1} \alpha_{2} \dotsm \alpha_{m} \) of positive
    integers whose sum is \( n \).  We consider the empty sequence
    with no terms to be the unique composition of 0.  We will
    sometimes write compositions as sums rather than as words, as in
    $\alpha_1 + \alpha_2 + \dotsb + \alpha_m$, though it must be kept
    in mind that the order of the terms still matters.  It is
    sometimes helpful to think of a composition of \( n \) as a
    sequence of \( n \) stones laid in a row, together with a grouping
    together of the stones in such a way that every stone belongs to a
    group, every group contains a stone, no stone belongs to two
    groups, and two stones belong to the same group only if every
    stone between them belongs to that group.

    Each term $\alpha_{i}$ in a composition \( \alpha \) is called a
    {\em part} of that composition.  A part equal to \( k \) is called
    a \( k \)-part.  A \( split \) in a composition is an integer
    that can be expressed as the sum of the first \( i \) parts of
    the composition for some nonnegative integer \( i \).  Thus, the
    composition $3 + 1 + 1 + 2$ has 5 splits: 0, 3, 4, 5, and 7.
    Using the imagery of stones, the groups of stones are the parts of
    \( \alpha \), and the splits of \( \alpha \) correspond to the
    spaces between groups (including the space to the left of the
    first group and to the right of the last group).  We use \(
    \rev{\alpha} \) to denote the composition comprising the parts of
    \( \alpha \) written in reverse order.  A {\em palindromic
    composition}, or a {\em palindrome}, is a composition for which \(
    \alpha = \rev{\alpha} \).  A {\em rise} (resp.  {\em drop}) is a
    part followed by a larger (resp.  smaller) part.  A {\em level} is
    a part followed by a part equal to itself.

    Frequencies of occurrences of \( k \)-parts, rises, drops, and
    levels in (palindromic) compositions, as well as in compositions
    with additional restrictions, have been studied (e.g.,
    see~\cite{CGH03} and~\cite{HM} and references therein).  Heubach
    and Mansour (\cite{HM}) give a multivariate generating function
    for joint distribution of parts, rises, levels, and drops in
    compositions and palindromes.  However, using the results from the
    literature related to the subject, it does not seem to be possible
    to answer a question like: how many levels are followed by rises
    among all compositions of \( n \)?  To consider a more general
    question, we introduce the notion of a {\em segmented pattern} in
    a composition.  A segmented pattern is a word $w = w_1 w_2 \cdots
    w_k$ in the alphabet of positive integers such that if \( b \) is
    a letter in \( w \) and \( a < b \), then \( a \) is a letter of
    \( w \).  In other words, the letters in \( w \) constitute an
    order ideal.  For example, $431242$ is a segmented pattern, while
    $41242$ is not.  We say that $w$ occurs in a composition \( \alpha
    = \alpha_{1} \alpha_{2} \dotsb \alpha_{m} \) if there is a subword
    $\alpha_{i} \alpha_{i+1} \dotsc \alpha_{i+k-1}$ of $\alpha$ that
    is order-isomorphic to $w$.  Thus, rises, drops, and levels are
    occurrences of the patterns 12, 21, and 11, respectively.  A level
    followed by a rise is an occurrence of the pattern 112.  The
    notion of segmented patterns in arbitrary words was studied in
    \cite{BM03}.

    More generally, we study occurrences of so-called {\em segmented
    partially ordered patterns} (SPOPs) in compositions.  A SPOP \( w
    \) is a word consisting of letters from a partially ordered
    alphabet \( \mathcal{A} \) such that the letters in \( w \)
    constitute an order ideal in \( \mathcal{A} \).  For instance, if
    we have a poset on three elements labeled by \( 1 \), $1'$, and
    $2'$ in which the only relation is $1' < 2'$, then the sequence
    31254 has two occurrences of $1 1' 2'$, namely $3 1 2$ and $1 2
    5$.  Given a SPOP \( w = w_{1} w_{2} \dotsm w_{m} \), we say that
    a segmented pattern \( v = v_{1} v_{2} \dotsm v_{m} \) is a {\em
    linear extension} of \( w \) if \( w_{i} < w_{j} \) implies that
    \( v_{i} < v_{j} \).  Thus the linear extensions of $11'2'$
    are $123$, $213$, and $312$.

    This paper is organized as follows.  In Section
    \ref{compositions}, we give our main results.  Theorem
    \ref{mainresult} gives a multivariate generating function for the
    number of occurrences of a given SPOP at a given split among
    compositions of $n$ with a given number of parts.  By specializing
    variables, we obtain a generating function for the number of
    occurrences of a given SPOP among all compositions of \( n \)
    (Corollary \ref{maincor}).  In Theorem \ref{numkpartsthm}, we
    enumerate the occurrences of \( k \)-parts at a given split in
    compositions of \( n \) with a given number of parts.  This
    generalizes a result in \cite{CCFM92}.  Our approach to this
    problem is to use a method which perhaps can be best described as
    ``encoding by restricted permutations."  The idea here is to
    encode a set of objects under consideration as a set of
    permutations satisfying certain restrictions.  Under appropriate
    encodings, this allows us to transfer the interesting statistics
    from our original set to the set of permutations, where they are
    easy to handle.  In Section \ref{palindromes}, we use restricted
    permutations to enumerate \( k \)-blocks with certain statistics
    in palindromic compositions, refining results in \cite{CGH03}.  In
    Section \ref{last} we provide short bijective encodings of binary
    bitonic sequences, binary strings without singletons, permutations
    avoiding $1\mn 3\mn 2\mn 4$ and having exactly one descent, and
    lines drawn through the points of intersections of $n$ straight
    lines in a plane.  Relations of these objects to certain
    restricted permutations were given in~\cite{BK}, but no bijections
    were provided.  We believe that these examples provide some
    evidence for the broad applicability of the method of encoding by
    restricted permutations.

    We use the following notations throughout the paper.  The set of
    nonnegative integers is denoted by \( \N \), and the set of
    positive integers is denoted by \( \Pos \).  Given \( m \leq n \in
    \N \), we write $[m,n]=\{m, m+1,\ldots, n\}$ and $[n] = [1,n]$.
    The permutations in this paper are written in one-line notation.
    Given a generating function \( G(t) \), we write \( [t^{n}] G(t)
    \) to denote the coefficient of \( t^n \) in \( G(t) \).  We use
    \( \comp(n) \) to denote the set of compositions of \( n \), and
    we write that \( \abs{\alpha} = n \) if \( \alpha \in \comp(n) \).
    Finally, let \( C(n,\ell) \) be the number of compositions of \( n
    \) with \( \ell \) parts.  It is well known and easy to verify
    that for a fixed non-negative integer \( \ell \), the generating
    function for \( C(n,\ell) \) is given by
    \begin{equation} \label{genfuncforcompswparts}
        \sum_{n=0}^{\infty} C(n,\ell) x^{n}
        =
        \frac{x^{\ell}}{(1-x)^{\ell}}.
    \end{equation}

    \section{Compositions} \label{compositions}

    Given a SPOP \( w = w_{1}w_{2} \dotsm w_{m} \) with \( m \) parts,
    let \( c_{w}(n,\ell,s) \) be the number of occurrences of \( w \)
    among compositions of \( n \) with \( \ell + m \) parts such that
    the sum of the parts preceding the occurrence is \( s \).  Let \(
    \Omega_{w}(x,y,z) \) be the generating function for \(
    c_{w}(n,\ell,s) \):
    \[
            \Omega_{w}(x,y,z)
            =
            \sum_{n,\ell,s \in \N} c_{w}(n,\ell,s) x^{n} y^{\ell} z^{s}.
    \]
    Our goal is to derive an explicit rational function for \(
    \Omega_{w}(x,y,z) \).

    Before proceeding, we define the following notation.  Given a
    segmented pattern \( v \) and \( n \in \N \), let \( P_{v}(n) \)
    denote the number of compositions of \( n \) that are order
    isomorphic to \( v \).  The generating function \(
    \mathcal{P}_{v}(x) \) for \( P_{v}(n) \) is not difficult to
    derive.  If \( j \) is the largest letter of \( v \), then \(
    P_{v}(n) \) is the number of integral solutions \( t_{1}, \dotsc,
    t_{j} \) to the system
    \begin{equation} \label{polytope}
        \mu_{1} t_{1} + \dotsb + \mu_{j} t_{j} = n,
        \qquad
        0 < t_{1} < \dotsb < t_{j},
    \end{equation}
    where \( \mu_{k} \) is the number of \( k \)'s in \( v \).  By
    expanding terms into geometric series, one can see that the number
    of integral solutions to \eqref{polytope} is the coefficient of \(
    x^{n} \) in
    \begin{equation} \label{genfunction}
        \mathcal{P}_{v}(x)
        =
        \prod_{k=1}^{j} \frac{x^{m_{k}}}{1 - x^{m_{k}}},
    \end{equation}
    where \( m_{k} = \mu_{j-k+1} + \dotsb + \mu_{j} \) for
    \( 1 \leq k \leq j \).

    \begin{thm} \label{mainresult}
        Let \( w \) be a SPOP.  Then
        \begin{equation} \label{genfunc}
            \Omega_{w}(x,y,z)
            =
            \sum_{v}
            \frac{(1-x)(1-xz)\mathcal{P}_{v}(x)}{(1-x-xy)(1-xz-xyz)}
        \end{equation}
        where the sum is over all linear extensions \( v \) of \( w
        \).
    \end{thm}

    \begin{proof}
        We begin by computing \( \Omega_{v}(x,y,z) \) when \( v \) is
        a segmented pattern.  We think of an occurrence of \( v \) as
        the triple of compositions \( (\alpha, \beta, \gamma) \) such
        that \( \alpha \) comprises the parts to the left of the
        occurrence, \( \beta \) comprises the parts in the occurrence,
        and \( \gamma \) comprises the parts to the right of the
        occurrence.  Hence, for given \( n, \ell, s \in \N \), \(
        c_{v}(n,\ell,s) \) is the number of triples \( (\alpha, \beta,
        \gamma) \) of compositions such that \( \abs{\alpha} +
        \abs{\beta} + \abs{\gamma} = n \), \( \abs{\alpha} = s \), \(
        \beta \) is order isomorphic to \( v \), and \( \alpha \) and
        \( \gamma \) together have \( \ell \) parts.  Thus we have
        that
        \[
            c_{v}(n,\ell,s)
            =
            \sum_{\substack{0 \leq j \leq \ell \\ 0 \leq k \leq n}}
            C(s,j) P_{v}(k) C(n-s-k,\ell-j).
        \]
        Using this equality, together with the generating function
        \eqref{genfuncforcompswparts} for \( C(n,\ell) \), we can
        factor \( \Omega_{v}(x,y,z) \) into a product of \(
        \mathcal{P}_{v}(x) \) and two geometric series:
        \begin{align*}
            \Omega_{v}(x,y,z)
            & = \sum_{n,\ell,s \in \N}
                c_{v}(n,\ell,s) x^{n} y^{\ell} z^{s} \\
            & = \mathcal{P}_{v}(x)
                \Bigg(
                    \sum_{n,\ell \in \N}
                    C(n,\ell) x^{n} y^{\ell}
                \Bigg)
                \Bigg(
                    \sum_{s,\ell \in \N}
                    C(s,\ell) (xz)^{s} y^{\ell}
                \Bigg) \\
            & = \mathcal{P}_{v}(x)
                \left(
                    \sum_{l \in \N}
                    \frac{x^{\ell}}{(1-x)^{\ell}} y^{\ell}
                \right)
                \left(
                    \sum_{l \in \N}
                    \frac{(xz)^{\ell}}{(1-xz)^{\ell}} y^{\ell}
                \right) \\
            & = \mathcal{P}_{v}(x)
                \frac{(1-x)(1-xz)}{(1-x-xy)(1-xz-xyz)}.
        \end{align*}
        Finally, note that if \( w \) is a SPOP, then \( c_{w}(x,y,z)
        = \sum_{v} c_{v}(x,y,z) \), where the sum is over all linear
        extensions \( v \) of \( w \).  Thus, \( \Omega_{w}(x,y,z) =
        \sum_{v} \Omega_{v}(x,y,z) \), and the theorem follows.
    \end{proof}

    Setting \( y=z=1 \) in equation \eqref{genfunc} yields the
    following.

    \begin{cor} \label{maincor}
        Given a segmented pattern \( w \), the number of occurrences of
        \( w \) among compositions of \( n \) is equal to
        \[
            [x^{n}]\Omega_{w}(x,1,1)
            =
            [x^{n}]
            \sum_{v} \frac{(1-x)^{2} \mathcal{P}_{w}(x)}{(1-2x)^{2}},
        \]
        where the sum is over all linear extensions \( v \) of \( w
        \).
    \end{cor}

    \begin{exm}
        We compute the number of occurrences of \( m \) levels
        immediately followed by a rise.  This is an occurrence of the
        segmented pattern \( w = \underbrace{1 \dotsm 1}_{m+1}2 \).
        The content vector of \( w \) is \( \mu = (m+1, 1) \), so we
        have
        \[
            \mathcal{P}_{w}(x)
            =
            \frac{x^{m+3}}{(1-x)(1-x^{m+2})}.
        \]
        Hence, the number of occurrences of \( w \) among all
        compositions of \( n \) is
        \[
            [x^{n}]\Omega_{w}(x,1,1)
            =
            [x^{n}]\frac{(1-x)^{2} x^{m+3}}{(1-2x)^{2} (1-x) (1-x^{m+2})}
        \]
        For fixed \( m \), it is routine to expand the rational
        function above into partial fractions to obtain a closed form
        expression for \( [x^{n}]\Omega_{w}(x,1,1) \).
    \end{exm}

    We now give an enumerative result that describes the number of \(
    k \)-parts located at a given split among compositions of \( n \)
    with a given number of parts.  Theorem \ref{numkpartsthm} below is
    our first example of encoding with restricted
    permutations.\footnote{In fact it is possible to use this result
    to prove Theorem \ref{mainresult}, though this approach requires
    several pages of tedious calculation, and is omitted in favor of
    the short and self-contained proof given above.} For \( n, k,
    \ell, s \in \N\), define \( f(n, k, \ell, s) \) to be the number
    of \( k \)-parts occurring among compositions of \( n \) with \(
    \ell + 1 \) parts such that the sum of the parts preceding the
    $k$-part is \( s \).  It immediately follows that \( f(n, k, \ell,
    s) = 0 \) if either \( n = 0 \) or \( k = 0 \).  The case when \(
    n = k > 0 \) is also clear: \( f(n, n, \ell, s) = 1 \) if \( \ell
    = s = 0 \), and \( f(n, n, \ell, s) = 0 \) otherwise.  The
    following theorem gives the value of \( f(n, k, \ell, s) \) in all
    remaining cases.

    \begin{thm} \label{numkpartsthm}
        If \( n \in \Pos \) and \( k \in [n-1] \), then
        \begin{equation} \label{numkparts}
            f(n, k, \ell, s)
            =
            \begin{cases}
                \binom{ n - k - 1}{\ell - 1}, & \text{if \( s \in \{0, n-k\} \) and \( \ell \in [n-k] \)}, \\
                \binom{ n - k - 2}{\ell - 2}, & \text{if \( s \in [n-k-1] \) and \( \ell \in [2, n-k] \)}, \\
                0,                            & \text{otherwise}.
            \end{cases}
        \end{equation}
    \end{thm}

    \begin{proof}
        We give a bijection between the \( k \)-parts that we are
        enumerating and a particular set of restricted permutations.
        Let \( S \) be the set of permutations of the quotient group
        \( Z = \Z/(n-k+1)\Z \) of the form \( s w_{1} w_{2} \dotsb
        w_{n-k} \), where
        \[
            w_{1} > w_{2} > \dotsb
            > w_{\ell}
            < w_{\ell+1} < \dotsb < w_{n-k}.
        \]
        and \( s + 1 \in \{w_{1}, \dotsc, w_{\ell}\} \) (where we've
        identified \( s \) with its canonical projection in \( Z \)).
        To see that \( \abs{S} \) is given by the right hand side of
        equation \eqref{numkparts}, observe that an element of \( S \)
        is uniquely specified by choosing which elements of \( Z \)
        will be in \( \{ w_{1}, \dotsc, w_{\ell} \} \) other than \( s
        \) (which cannot be in there) and \( s + 1 \) and \( \min(Z
        \backslash \{s\}) \) (which must be in there, but which are
        equal when \( s \in \{0, n-k\} \)).

        We now show that the elements of \( S \) are in bijective
        correspondence with the \( k \)-parts that we wish to
        enumerate.  First, we can think of such a \( k \)-part as an
        element of
        \[
            T
            =
            \{
                (\alpha, \beta)
                :
                \alpha \in \comp(s) ,\,
                \beta \in \comp(n-k-s) ,\,
                \text{
                    \( \alpha \) and \( \beta \) together have \( \ell
                    \) parts
                }
            \}.
        \]
        To be precise, a \( k \)-part in a composition of \( n \)
        corresponds to the ordered pair \( (\alpha, \beta) \) of
        compositions such that \( \alpha \) comprises the parts to the
        left of the chosen \( k \)-part and \( \beta \) comprises the
        parts to the right of the chosen \( k \)-part.

        We now give a bijection \( T \leftrightarrow S \).  For an
        explicit example of the bijection we are about to describe,
        see Example \ref{bijexample}.  Given \( (\alpha, \beta) \in
        T\), we produce a permutation in \( S \) as follows.
        Concatenate the compositions \( \alpha \) and \( \beta \),
        producing a composition \( \gamma \in \comp(n-k) \) with \(
        \ell \) parts.  Let \( \overline{w}_{\ell} = 0 \) and let
        \[
            \overline{w}_{\ell-i}
            =
            \sum_{j=1}^{i} \gamma_{j},
            \quad
            \text{for \( 1 \leq i \leq \ell-1 \)},
        \]
        For \( 1 \leq i \leq \ell \) let
        \[
            w_{i}
            =
            \begin{cases}
                \overline{w}_{i},     & \text{if \( \overline{w}_{i} < s \)},\\
                \overline{w}_{i} + 1, & \text{if \( \overline{w}_{i} \geq s\)}.
            \end{cases}
        \]
        Finally, let \( w_{\ell+1}, \dotsc, w_{n-k} \) be the elements
        of \( [0, n-k] \backslash \{s, w_{1}, \dotsc, w_{\ell}\} \)
        written in increasing order (where we've identified \( [0,n-k]
        \) with its canonical projection in \( Z \)).

        It is easy to show that this map yields an element of \( S \).
        We show that the map is a bijection by giving its inverse.
        Given an element of \( S \), one may produce an element of \(
        T \) by letting
        \[
            \overline{w}_{i}
            =
            \begin{cases}
                w_{i},     & \text{if \( w_{i} < s \)},\\
                w_{i} - 1, & \text{if \( w_{i} > s\)},
            \end{cases}
            \qquad
            \text{for \( 1 \leq i \leq \ell \)},
        \]
        letting \( \gamma_{i} = \overline{w}_{\ell-i} -
        \overline{w}_{\ell-i+1} \) for \( 1 \leq i \leq \ell \), and
        letting \( \gamma = \gamma_{1} \dotsm \gamma_{\ell} \).
        Because of the requirement that \( s+1 \in \{w_{1}, \dotsc,
        w_{\ell}\} \), it follows that for some \( i \), \(
        \sum_{j=1}^{i} \gamma_{i} = s \).  Let \( \alpha = \gamma_{1}
        \dotsm \gamma_{i} \) and let \( \beta = \gamma_{i+1} \dotsm
        \gamma_{\ell} \).  Then we have that \( (\alpha, \beta) \in T
        \).
    \end{proof}

    \begin{exm}\label{bijexample}
        We choose as our \( k \)-part the 6 in the composition \( 3 \,
        1 \, 6 \, 2 \).  Then we have \( n = 12 \), \( k = 6 \), \( l
        = 3 \), and \( s = 4 \).  The claim is that this corresponds
        to a permutation of the elements in \( \Z/7\Z \).

        Applying the maps from the theorem to our chosen \( k \)-part
        yields \( \alpha = 3 \, 1 \), and \( \beta = 2 \).  Thus we
        have \( \gamma = 3 \, 1 \, 2 \).  Computing the values of \(
        \overline{w}_{i} \) yields \( \overline{w}_3 = 0 \), \(
        \overline{w}_2 = \gamma_1 = 3 \), and \( \overline{w}_1 =
        \gamma_1 + \gamma_2 = 4 \).  Observing that \( \overline{w}_1
        \geq s = 4 \), while \( \overline{w}_2, \overline{w}_3 < s \),
        we compute the \( w_{i} \)'s as follows
        \begin{align*}
            w_1 & = \overline{w}_1 + 1 = 5, \\
            w_2 & = \overline{w}_2 = 3, \\
            w_3 & = \overline{w}_3 = 0.
        \end{align*}
        Finally, we let \( w_4 \, w_5 \, w_6 \) be the elements of
        \begin{align*}
            \{0, . . ., 6\} \backslash \{s, w_1, w_2, w_3\}
            & = \{0, . . . , 6\} \backslash \{4, 5, 3, 0\} \\
            & = \{1, 2, 6\}
        \end{align*}
        written in increasing order.  Therefore, the word
        corresponding to our original \( k \)-part is
        \[
            s \, w_1 \,  w_2 \dotsm w_6 = 4 \, 5 \, 3 \, 0 \, 1 \, 2 \, 6.
        \]
    \end{exm}

    As a corollary to Theorem \ref{numkpartsthm}, we derive a result
    that appeared in \cite{CCFM92}.

    \begin{cor} \label{numkpartscor}
        Given \( n \in \N \) and \( k \in [n-1] \), the number of \( k
        \)-parts among all compositions of \( n \) is \( 2^{n-k-2} (n
        - k + 3) \).
    \end{cor}

    \begin{proof}
        The result follows from using equation \eqref{numkparts} to
        compute
        \begin{align*}
            \sum_{\substack{\ell \in [n-k] \\ s \in [0, n-k]}} f(n,k,\ell,s)
            & =
            2 \sum_{\ell=1}^{n-k} \binom{n-k-1}{\ell-1}
            +
            (n - k - 1) \sum_{\ell=2}^{n-k} \binom{n-k-2}{\ell-2} \\
            & =
            2^{n-k} + 2^{n-k-2} (n-k-1) \\
            & =
            2^{n-k-2} (n - k + 3)
        \end{align*}
    \end{proof}

    \section{Palindromic Compositions}\label{palindromes}

    We provide two alternative (nonequivalent) encodings by restricted
    permutations of $k$-parts in palindromes of $N$, when $N$ and $k$
    have different parity.  We give these encodings explicitly in the
    case of even palindromes of $N=2(n-1)$ and odd $k$-parts.  These
    encodings provide bijective proofs of the known result that the
    number of \( k \)-parts in palindromic compositions of $2(n-1)$ is
    $(n-k+1)2^{n-k-1}$ when \( k \) is odd (see \cite{CGH03}).  Such \( k
    \)-parts will be encoded as permutations \( w_{1}
    w_{2}\, \dotsm w_{n-k+1} \) of \( [n-k+1] \) such that, for some
    \( \ell \in \{2, \dotsc, n\} \), $w_{2} > w_{3} > \dotsb >
    w_{\ell} < w_{\ell+1} < \dotsb < w_{n-k+1}$.  In either encoding,
    the case of odd palindromes $N=2n-1$ and even $k$-parts can be
    obtained using similar ideas.

    \subsection{First encoding}

    First, observe that a permutation \( w_{1} w_{2} \dotsm w_{n-k+1}
    \) of \( [n-k+1] \) with
    \[
        w_{2} > w_{3} > \dotsb > w_{\ell} < w_{\ell+1} < \dotsb < w_{n-k+1}.
    \]
    corresponds to an ordered pair \( (w_{1},\alpha) \) with \( w_{1}
    \in [n-k+1] \) and \( \alpha = \alpha_{1} \dotsm \alpha_{\ell-1} \in
    \comp(n-k) \) as follows.  If \(\ell = 2 \), let \( \alpha_{1} = n-1
    \).  Otherwise, let
    \[
        \overline{w}_{i}
        =
        \begin{cases}
            w_{i} - 1, & \text{if \( w_{i} < w_{1} \)}, \\
            w_{i} - 2, & \text{if \( w_{1} < w_{i} \)},
        \end{cases}
        \qquad
        \text{for \( 2 \leq i \leq l-1 \)},
    \]
    and put
    \begin{align*}
        \alpha_{1}      & = \overline{w}_{\ell-1}, \\
        \alpha_{i}      & = \overline{w}_{\ell-i} - \overline{w}_{\ell-i+1},
                            \quad
                            \text{for \( 2 \leq i \leq l-2 \)}, \\
        \alpha_{\ell-1} & = n - k - \overline{w}_{2}.
    \end{align*}
    We now explicitly describe the correspondence between pairs \(
    (w_{1}, \alpha) \) and odd \( k \)'s in palindromic compositions
    of \( 2(n-1) \).  It may be helpful to use the imagery of stones
    discussed in Section \ref{intro}.  In this context, \(
    w_{1} \) can be thought of as distinguishing a {\em gap} in the
    sequence of stones, where the gaps are the spaces between
    any two adjacent stones (whether they belong to the same group or
    not), as well the space before the first stone and after the last
    stone.  Hence, a sequence of \( n-k\) stones has \( n-k+1 \) gaps,
    which are indexed with the set \( [n-k+1] \).

    \bigskip

    \noindent \textbf{Case I}: The cases in which \( w_{1} \in \{1,
    n-k+1\} \) correspond to the \( k \)'s that are either the
    left-most or right-most terms in the compositions containing them.
    In particular, \( (1, \alpha) \) corresponds to the left-most \( k
    \) in the composition
    \[
        k
        +
        \sum_{i=1}^{\ell-2} \alpha_{i}
        +
        2(\alpha_{\ell-1} - 1)
        +
        \rev{\sum_{i=1}^{\ell-2} \alpha_{i}}
        +
        k,
    \]
    while \( (n-k+1, \alpha) \) corresponds to the right-most \( k \).

    \bigskip

    \noindent {\bfseries Case II}: The cases in which \( 2 \leq w_{1}
    \leq n-k \) and \( w_{1} \) is a split in \( \alpha \) correspond to
    the \( k \)'s which are on the left-hand side of the palindromic
    compositions containing them, but which are not the left-most
    terms.  In these cases, \( (w_{1}, \alpha) \) corresponds to the
    indicated \( k \) on the left-hand side of the palindromic
    composition
    \[
        \sum_{i=1}^{j} \alpha_{i}
        +
        k
        +
        \sum_{i=j+1}^{\ell-2} \alpha_{i}
        +
        2(\alpha_{\ell-1} - 1)
        +
        \rev{\sum_{i=1}^{j} \alpha_{i}
        +
        k
        +
        \sum_{i=j+1}^{\ell-2} \alpha_{i}},
    \]
    where we have the identity \( \sum_{i=1}^{j} \alpha_{i} = w_{1} -
    1 \).

    \bigskip

    \noindent {\bfseries Case III}: The cases in which \( 2 \leq w_{1}
    \leq n-k \) and \( w_{1} \) is not a split in \( \alpha \)
    correspond to the \( k \)'s which are on the right-hand side of
    the palindromic compositions containing them, but which are not
    the right-most terms.  These cases break into two subordinate
    cases:

    \bigskip

    \noindent {\bfseries Case IIIA}: Within Case III, those \( (w_{1},
    \alpha) \) in which \( w_{1} \) is a gap in the last term of \(
    \alpha \) correspond to \( k \)'s that are immediately to the
    right of the center term of the palindromic compositions
    containing them.  In particular, such a \( (w_{1}, \alpha) \)
    corresponds to the indicated \( k \) on the right-hand side of the
    palindromic composition
    \[
        \sum_{i=1}^{\ell-2} \alpha_{i}
        +
        \alpha_{\ell-1}'
        +
        k
        +
        2(\alpha_{\ell-1}'' - 1)
        +
        k
        +
        \alpha_{\ell-1}'
        +
        \rev{\sum_{i=1}^{\ell-2} \alpha_{i}},
    \]
    where we use the identities \( \alpha_{\ell-1}' +
    \alpha_{\ell-1}'' = \alpha_{\ell-1} \) and \( \sum_{i=1}^{\ell-2}
    \alpha_{i} + \alpha_{\ell-1}' = w_{1} - 1 \).

    \bigskip

    \noindent {\bfseries Case IIIB}: On the other hand, if \( w_{1} \)
    is not a gap in the last term of \( \alpha \), then \( (w_{1},
    \alpha) \) corresponds to the indicated \( k \) on the right-hand
    side of the palindromic composition
    \begin{multline*}
        \sum_{i=1}^{j-1} \alpha_{i}
        +
        \alpha_{j}'
        +
        k
        +
        \alpha_{j}''
        +
        \sum_{i=j+1}^{\ell-2} \alpha_{i}
        +
        2(\alpha_{\ell-1} - 1) \\
        +
        \rev{\sum_{i=1}^{j-1} \alpha_{i}
        +
        \alpha_{j}'
        +
        k
        +
        \alpha_{j}''
        +
        \sum_{i=j+1}^{\ell-2} \alpha_{i}},
    \end{multline*}
    where we use the identities \( \alpha_{j}' + \alpha_{j}'' =
    \alpha_{j} \) and \( \sum_{i=1}^{j-1} \alpha_{i} + \alpha_{j}' =
    w_{1} - 1 \).  (To deduce the value of \( j \) from a given
    composition \( \alpha \), we will also need to use the inequality
    \( 0 < \alpha_{j}' < \alpha_{j}\).)

    \subsection{Second Encoding}

    Clearly, a palindrome of $2(n-1)$ has either an odd number of
    parts with an even part in the center or an even number of parts
    and no central part.  To make all palindromes to be of odd length,
    we create a central part ``0" for palindromes with an even number
    of parts.

    We present an algorithm to produce a permutation given an
    underlined $k$-part in a palindrome $P$.  We only consider the
    case when the chosen part is to the left of the center in $P$; for
    a part from the right-hand side, we proceed with the part
    symmetric to it, and we switch 1 and 2 in the obtained
    permutation.  In the bijection below, a part is to the left of the
    center if and only if in the corresponding permutation, 1 precedes
    2.  In general, we find the permutation corresponding to
    $\underline{k}$ by inserting the numbers $n-k+1$, $n-k$, $n-k-1$,
    and so on in decreasing order, into initially empty slots $w_1,
    w_2\ldots, w_{n-k+1}$.

    Suppose $P=C\underline{k}DxDkC$ where $x=2t$ for $t\geq 0$.
    \begin{enumerate}
        \item
        If $D$ is empty and $x=0$, we set $w_1=(n-k+1)$ and proceed
        with (2) below.  Otherwise, we set
        $$w_{n-k-t+2}w_{n-k-t+3}\dotsm
        w_{n-k+1}=(n-k-t+2)(n-k-t+3)\dotsm (n-k+1)$$ and
        $w_2=(n-k-t+1)$ (in particular, if $t=0$ we only set
        $w_2=(n-k+1)$).  We read the parts in $D$ from right to left
        and fill in the slots $w_3,w_4,\ldots$ by placing $n-k-t$,
        then $n-k-t-1$, and so on: if a current part is $a$, then we
        place $a-1$ of the largest unplaced numbers to the right in
        increasing order, and we place the largest number out of the
        remaining numbers to the left.  We then proceed with the part
        next to $a$ from the left.  The only exception is the part
        immediately to the right of $\underline{k}$.  In this case, we
        place $a-1$ of the largest unplaced numbers to the right in
        increasing order, and then we set $w_1$ be the largest of yet
        unplaced numbers.  If we get $w_2=2$, set $w_1=1$ and place 2
        in the only one remaining slot.  Continue with step (2).

        \item
        If $C$ is empty or if $C=1$, place the unplaced numbers in
        increasing order into the empty slots.  Otherwise, suppose
        $C=a_1 a_2 \cdots a_k$.  Then we consider the binary vector
        $0^{a_1-1}10^{a_2-1}1\cdots 0^{a_k-1}$ (each block of 0's but
        the last one is followed by a 1).  We read this binary vector
        from right to left and whenever we meet a 0, we place the
        largest unplaced number into the leftmost available slot;
        otherwise, we place this number into the rightmost available
        slot.  If this procedure can no longer be continued, and 1 or
        2 have not yet been placed, place them so that 1 precedes 2.
    \end{enumerate}

    We provide some examples.  Suppose we are interested in
    $\underline{1}$ in the following palindrome of 16:
    $2\underline{1}2141212$.  The steps of our recursive bijection are
    as follows:
    $$
        *7*****89
        \rightarrow
        *76****89
        \rightarrow
        476***589
        \rightarrow
        4763**589
        \rightarrow
        476312589.
    $$
    As further examples, one can check that the underlined 5's in
    $\underline{5}115$, $15\underline{5}1$, and $\underline{5}25$
    correspond to 132, 321, and 123 respectively.

    The inverse of this algorithm is easy to find.  In particular, if
    $w_1=1$ (resp.  $w_1=2$) then the corresponding $k$-part is the
    leftmost (resp.  rightmost) one in a composition.

    \section{Additional Encodings with Restricted Permutations} \label{last}

    We now provide some additional examples of encodings of
    combinatorial objects by restricted permutations to demonstrate
    various approaches to bijective enumeration.  But first, some
    definitions.

    A sequence $a_1, a_2, \ldots, a_n$ is \emph{bitonic} if for some
    $h$, $1\leq h\leq n$, we have that $a_1\le a_2\le\cdots\le a_h\ge
    a_{h+1}\ge \cdots \ge a_{n-1}\ge a_n$ or $a_1\ge a_2\ge\cdots\ge
    a_h\le a_{h+1}\le \cdots \le a_{n-1}\le a_n$.  A binary string $x$
    is said to be {\em without singletons} if the words 010 and 101
    are not factors of $x$.

    Let $S_1$ (resp.  $S_2$) be the set of $(n+2)$-permutations $w_{1}
    w_{2}\, \dotsm w_{n+2}$ such that, $w_{1}w_{2}=(n+1)(n+2)$ or
    $w_{1}w_{2}=(n+2)(n+1)$, and $w_{3}w_{4}\dotsm w_{n+2}$ avoids
    simultaneously the patterns $1\mn 2\mn 3$ and $2\mn 3\mn 1$ (resp.
    $1\mn 2\mn 3$, $1\mn 3\mn 2$, and $2\mn 1\mn 3$).  According
    to~\cite{BK}, $\abs{S_1} = n^2-n+2$ and $\abs{S_2}=2F_n$, where
    $F_n$ is the $n$-th Fibonacci number with $F_0=F_1=1$.

    Let $S_3$ be the set of $(n+3)$-permutations $w_{1} w_{2}\, \dotsm
    w_{n+3}$ such that, $w_{1}<w_{2}<w_{3}$ and $w_{4}w_{5}\dotsm
    w_{n+3}$ is in decreasing order.  Clearly, $\abs{S_3} =
    \binom{n+3}{3}$.

    Let $S_4$ be the set of $n$-permutations $w_{1} w_{2}\, \dotsm
    w_{n}$ such that, $w_1$ is the largest letter among the four
    leftmost letters, $w_{3}<w_{4}$ and $w_{5}w_{6}\dotsm w_{n}$ is in
    decreasing order.  One can see that $\abs{S_4} = 3\binom{n}{4}$.

    \begin{bijection}
        The elements of $S_1$ are in one-to-one correspondence with
        binary bitonic sequences of length $n-1$.
    \end{bijection}

    In order to avoid the restrictions, $w_{3}w_{4}\dotsm w_{n+2}$
    must be either of the form $i(i-1)\dotsm 1n(n-1)\dotsm (i+1)$ or
    of the form $$n(n-1)\dotsm (n-i+1)(j+1)j\dotsm 1(n-i)(n-i-1)\dotsm
    (j+2)$$ for some $i>0$ and $j\geq 0$.

    We describe our bijection in the case $w_{1}w_{2}=(n+1)(n+2)$.  We
    then use the same bijection for $w_{1}w_{2}=(n+2)(n+1)$ and
    replace 0's by 1's and 1's by 0's in the corresponding sequences.

    To the permutation $(n+1)(n+2)i(i-1)\dotsm 1n(n-1)\dotsm (i+1)$
    there corresponds the bitonic sequence $01^{i}0^{n-i-2}$ where
    $i>0$; to the permutation
    $$
        (n+1)(n+2)n(n-1) \dotsm (n-i+1)(j+1)j \dotsm 1(n-i)(n-i-1) \dotsm (j+2)
    $$
    there corresponds the sequence $00^{i}1^{n-i-j-2}0^{j}$ where
    $i>0$ and $j\geq 0$.  Clearly, our map involves all the binary
    bitonic sequences starting from 0 exactly once and the reverse to
    this map is easy to see.  Together with the case
    $w_{1}w_{2}=(n+2)(n+1)$ we have a bijection.

    \begin{bijection}
        The elements of $S_2$ are in one-to-one correspondence with
        binary strings of length $n+2$ without singletons.
    \end{bijection}

    Clearly, any string under consideration ends with either 00 or
    with 11.  We match the strings ending with 00 with the
    permutations beginning with $w_1 w_2=(n+1)(n+2)$.  It will suffice
    to consider this case.  The remaining cases are handled by
    replacing 0's by 1's and 1's by 0's, proceeding with the first
    case, and then replacing \( (n+1)(n+2) \) with \( (n+2)(n+1) \) in
    the resulting permutation.

    We begin with a procedure for construction permutations
    $w_{3}w_{4}\dotsm w_{n+2}$ that avoid the restricted patterns.
    Insert the numbers $1,2,\ldots,n$, in that order, into $n$ slots
    corresponding to the letters $w_i$, $3\leq i\leq n+2$.  Note that
    we must either set $w_{n+2}=1$ or set $w_{n+1}w_{n+2}=12$, since
    otherwise we get an occurrence of a prohibited pattern.  We
    proceed by induction.  If the rightmost $i$ slots have been filled
    and \( w_{n-i} \) is empty, then we only have two choices: either
    set $w_{n-i}=i+1$ or set $w_{n-i-1}w_{n-i}=(i+1)(i+2)$.  A
    permutation $w_{3}w_{4}\dotsm w_{n+2}$ that avoids the restricted
    patterns may be thought of as a tiling of a \( 1 \times n \) board
    by monominos and dominos.

    Now, given such a tiling, we construct a binary string \( b_{1}
    b_{2} \dotsm b_{n} 00 \) corresponding to that tiling.  Read the
    tiling from right to left.  If the leftmost tile is a monomino,
    set \( b_{n}=0 \).  Otherwise, set \( b_{n-1}b_{n} = 11 \).  In
    general, if the last digit placed in the binary string was \(
    b_{i} = x \in \{0,1\} \), and the next unread tile is a monomino,
    read this tile and set \( b_{i-1} = x \).  Otherwise, if the next
    unread tile is a domino, read this domino and set \(
    b_{i-2}b_{i-1} = \bar{x}\bar{x} \), where \( \bar{x} \) is the
    binary complement of \( x \).  In this way, we avoid the
    possibility of creating singletons.

    This process is reversible: we read a binary word without
    singletons from right to left while tiling a \( 1 \times n \)
    board with monominos and dominos.  Whenever we meet
    $\bar{x}\bar{x}$ after passing $x$ in the binary string, we place
    a domino on the board.  Otherwise, we place a monomino.  The
    resulting tiling defines the corresponding permutation according
    to the construction described above.  Performing all of these
    correspondences yields the desired bijection.  For example, if \(
    w_{1} w_{2} \dotsm w_{9} = 896753412 \), then we produce \( b_{1}
    b_{2} \dotsm b_{9} = 110001100 \).

    \begin{bijection}
        The elements of $S_3$ are in one-to-one correspondence with
        $(n+2)$-permutations avoiding $1\mn 3\mn 2\mn 4$ and having
        exactly one descent (a descent is an $i$ such that
        $w_i>w_{i+1}$).
    \end{bijection}

\begin{figure}
        \begin{center}
            \includegraphics[width=180pt]{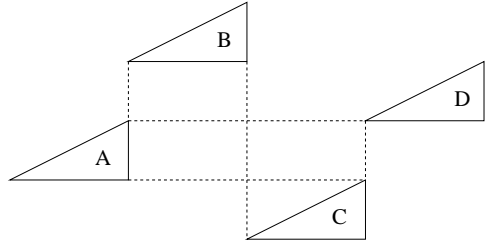}
            \caption{The structure of permutations avoiding $1\mn 3\mn
            2\mn 4$ and having exactly one descent.}
            \label{figure}
        \end{center}
    \end{figure}

    Any $(n+2)$-permutation avoiding $1\mn 3\mn 2\mn 4$ and having
    exactly one descent has the structure $ABCD$, where
    $A=(i_1+1)(i_1+2)\dotsm i_2$, $B=(i_3+1)(i_3+2)\dotsm (n+2)$,
    $C=12\dotsm i_1$, $D=(i_2+1)(i_2+2)\dotsm i_3$ (see
    Figure~\ref{figure}), and one of the following four mutually
    exclusive possibilities occurs:
    \begin{enumerate}
        \item
        none of $A$, $B$, $C$, and $D$ is empty: there are
        $\binom{n+1}{3}$ such permutations, given by the number of
        ways to choose the least elements in $A$, $B$, and $D$ (we
        know that 1 belongs to $C$);

        \item
        $C$ is empty: there are $\binom{n+1}{2}$ such permutations,
        since 1 belongs to $A$ and we choose the least elements in $B$
        and $D$;

        \item
        $B$ is empty: there are $\binom{n+1}{2}$ such permutations,
        since 1 belongs to $C$ and we choose the least elements in $A$
        and $D$;

        \item
        $A$ and $C$ are empty: there are $(n+1)$ such permutations,
        since $1$ is in $D$ and we need to choose the length of $D$
        ($B$ is not empty).
    \end{enumerate}

    Note that summing over all the cases gives us exactly
    $\binom{n+3}{3}$ permutations.  Once the permutations in avoiding
    \( 1 \mn 3 \mn 2 \mn 4 \) with exactly one descent have been
    partitioned into the four cases above, it is easy to find
    bijections in each case with permutations in \( S_{3} \) as
    follows.
    \begin{enumerate}
        \item
        $abc(n+3)(n+2)w_6w_7\dotsm w_{n+3}$, where $a<b<c$ and
        $w_6w_7\dotsm w_{n+3}$ is decreasing.  Choosing $a$, $b$, and
        $c$ corresponds to choosing $i_1$, $i_2$, and $i_3$;

        \item
        $ab(n+3)(n+2)w_5w_6\dotsm w_{n+3}$, where $a<b$ and
        $w_5w_6\dotsm w_{n+3}$ is decreasing.  Choosing $a$ and $b$
        corresponds to choosing $i_2$ and $i_3$;

        \item
        $ab(n+2)(n+3)w_5w_6\dotsm w_{n+3}$, where $a<b$ and
        $w_5w_6\dotsm w_{n+3}$ is decreasing.  Choosing $a$ and $b$
        corresponds to choosing $i_1$ and $i_2$;

        \item
        $a(n+2)(n+3)w_4w_5\dotsm w_{n+3}$, where $a<b$ and
        $w_4w_5\dotsm w_{n+3}$ is decreasing.  The length of $D$
        corresponds to $a$.
    \end{enumerate}
    Since these cases provide a partition of permutations in \( S_{3}
    \), the bijection is complete.

    \begin{bijection}
        The elements of $S_4$ are in one-to-one correspondence with
        the set of all lines drawn through the points of intersections
        of $n$ straight lines in a plane, no two of which are
        parallel, and no three of which are concurrent (we assume here
        that each of such lines goes through exactly two points of
        intersections).
    \end{bijection}

    If we label the lines by $1,2,\ldots,n$ then each intersection
    point can be represented by a pair of numbers $(x,y)$
    corresponding to the intersecting lines.  Now any line from the
    set of ``new" lines can be described by a pair $((x,y),(z,v))$
    where all of $x$, $y$, $z$, and $v$ are different.  Assuming that
    $x<y$, $z<v$, and $y<v$ we construct the corresponding permutation
    $vzxyw_{5}w_{6}\dotsm w_{n}$ where $w_{5}w_{6}\dotsm w_{n}$ is
    decreasing.  Clearly this map is a bijection.

\end{document}